\documentclass[11pt]{amsart}
\usepackage{titletoc}
\usepackage{setspace}
\usepackage[margin=1in]{geometry}

\numberwithin{equation}{section}
\usepackage{amssymb,latexsym,amsfonts}
\usepackage{pb-diagram}
\usepackage{graphicx}
\usepackage[pdftex]{hyperref}
\usepackage{verbatim}
\usepackage{color} 
\usepackage[matrix,arrow]{xy}
\usepackage{marginnote}
\usepackage{amsmath}

\newcommand{\scal}{\mathrm{scal}}

\theoremstyle{plain}
\newtheorem{theorem}{Theorem}[section]
\newtheorem{proposition}[theorem]{Proposition}
\newtheorem{lemma}[theorem]{Lemma}
\newtheorem{cor}[theorem]{Corollary}
\newtheorem{corollary}[theorem]{Corollary}

\newcommand{\sq}{\hbox{\rlap{$\sqcap$}$\sqcup$}}

\theoremstyle{definition}
\newtheorem{definition}[theorem]{Definition}

\DeclareMathOperator{\End}{End}

\newcommand{\Inv}{\mathcal{R}^{INV}}
\newcommand{\RicFlat}{\mathcal{RIC}_{=0}}

\newcommand{\integers}{{\mathbb Z}}
\newcommand{\complexs}{{\mathbb C}}

\newcommand{\iso}{\cong}
\newcommand{\tensor}{\otimes}

\begin{document}

\title{Non-negative versus positive scalar curvature}

\begin{abstract} 
  In this note, we look at the difference, or rather the absence of a
  difference, between the space of metrics of positive scalar curvature and
  metrics of non-negative scalar curvature. The main tool to analyze the
  former on a spin manifold is the spectral theory of the Dirac operator and
  refinements thereof. This can be used, for example, to distinguish between
  path components in the
  space of positive scalar curvature metrics. Despite the fact that
  non-negative scalar curvature a priori does not have the same spectral
  implications as positive scalar curvature, 
  we show that all invariants based on the Dirac operator extend over the
  bigger space. Under mild conditions we show that the inclusion of the
  space of metrics of positive scalar curvature into that of non-negative
  scalar curvature is a weak homotopy equivalence.
\end{abstract}

\author{Thomas Schick}
\address{
Mathematisches Institut\\
Universit\"at G\"ottingen \\
G\"ottingen \\
Germany}
\email{thomas.schick@math.uni-goettingen.de}

\author{David J. Wraith}
\address{Department of Mathematics and Statistics\\
National University of Ireland Maynooth\\
Maynooth\\
Ireland}
\email{david.wraith@mu.ie}

\date{\today}

\maketitle

\renewcommand{\baselinestretch}{1.3}\normalsize

\section{Introduction}\label{Intro}

The study of the topology of spaces and
moduli spaces of Riemannian metrics satisfying some form of curvature
condition on a fixed manifold has for many years been an important research
subject. Such curvature conditions include positive 
scalar curvature, positive Ricci curvature, and non-negative 
sectional curvature. For some recent results concerning closed manifolds, see for
example \cite{BHSW}, \cite{BERW}, \cite{HSS}, \cite{CS}, \cite{CSS}, \cite{CM},
\cite{BB}, \cite{Wal1}, \cite{Wal2},
\cite{Wr}, \cite{DKT}, \cite{TWi}, and the book \cite{TW}.

In this paper all manifolds under consideration will be closed and connected
unless otherwise stated, and we will always assume that spaces of metrics are
equipped with the $C^\infty$-topology.

The principal theme in this paper is the comparison of (moduli) spaces of
non-negative scalar curvature metrics with (moduli) spaces of positive scalar
curvature metrics on closed spin manifolds $M$. In this
context the Ricci flat metrics play a special role, and 
with this in mind we make the following

\begin{definition}
  Let $\mathcal N$ denote the space of non-negative scalar curvature metrics
  on $M$. Similarly, let
  $\mathcal P$ denote the space of positive scalar curvature metrics on
  $M$. Denote by $\RicFlat$ the space of Ricci flat metrics, and
  set $\mathcal{P}^\sharp:=\mathcal{P}\cup \RicFlat$.

  We have the obvious inclusion relations
  \begin{equation*}
    \mathcal{P}\subset \mathcal{P}^\sharp =\mathcal{P}\cup \RicFlat
    \subset \mathcal{N}.
  \end{equation*}
\end{definition}

{We claim that if $\mathcal{P}=\emptyset$, then $\mathcal{P}^\sharp=\RicFlat=\mathcal{N}$. To see this we begin by recalling that the Trichotomy Theorem of Kazdan and Warner (\cite{KW1},\cite{KW2}, compare \cite{RS}) implies that if $M$ admits a non-negative scalar curvature metric for which the scalar curvature is not identically zero, then $M$ in fact admits a positive scalar curvature metric. Thus if $\mathcal{P}=\emptyset$ and $g \in \mathcal{N},$ we conclude that $\scal(g) \equiv 0.$ A classical result of Bourguignon (compare \cite[4.49]{Be}) now asserts that if the only metrics on $M$ with non-negative scalar curvature are scalar flat, then any scalar flat metric on $M$ must be Ricci flat. This establishes the claim.}

{Using the Ricci flow, we see that in a homotopy theoretic sense the above claim remains true in general:}

  \begin{theorem}\label{maintheorem-a}
 The inclusion $\mathcal{P}^\sharp\hookrightarrow \mathcal{N}$ is a weak homotopy equivalence.
  \end{theorem}

  \begin{corollary}\label{theo:main-no_Ricci}
  Let $M$ a closed spin manifold which does not admit a Ricci flat Riemannian
  metric. Then the inclusion $\mathcal{P}\hookrightarrow \mathcal{N}$ is a
  weak homotopy equivalence.
\end{corollary}

Note that it is rare that a manifold admits a Ricci flat metric. For example,
by \cite{CG} (see also Theorem \ref{structure}), the fundamental 
group of a Ricci flat manifold contains a free abelian subgroup (possibly
trivial) of finite index.

In view of Corollary \ref{theo:main-no_Ricci}, the interesting case for our
investigation now is the complementary case
where $\mathcal{P}\ne \mathcal{P}^\sharp$, i.e.~$\RicFlat\ne\emptyset$.

Most of the results to date concerning (moduli) spaces of positive scalar
curvature metrics are established using the index theory of Dirac
operators. We will present some of the relevant details concerning this in
Section \ref{basicindex}, however for now it suffices to note that
one of the key results which makes index theory such an important tool in this
context is the classical theorem of Schr\"odinger-Lichnerowicz. In order to
state this, let us first recall that if $(M,g)$ is a Riemannian spin manifold,
we can consider the spin Dirac operator $\mathcal{D}$ defined by Atiyah and Singer acting on
the space of sections of the spinor bundle over $M$. This operator depends on
the metric and on the spin structure. Sections which belong to
the kernel of $\mathcal{D}$ are called harmonic spinors. The basic case of the Schr\"odinger-Lichnerowicz
Theorem then states that a compact spin manifold with positve scalar curvature
admits no non-trivial harmonic spinors.

One can extend this result by generalizing the concepts of Dirac operator and
harmonic spinor by
twisting the spinor bundle (that is, forming the tensor product) with a flat
bundle $F$ over the same base, see for example \cite[pages 164-165]{LM}. This is an
important construction, and is frequently used in the case where the flat bundle is
not a vector bundle, but a bundle of modules over an auxiliary $C^*$-algebra
$A$. For us, the most relevant $C^*$-algebra is the maximal (real) $C^*$-algebra of the
fundamental group of $M$, called $C^*\pi_1(M)$. In some sense there is a universal
case of twisting with a flat bundle, and this involves the so-called
Mishchenko line bundle $L_M$ over $M.$  This is
a bundle with fibre a free rank one 
module over $C^*\pi_1(M)$, which comes equipped with a canonical flat connection. We will recall the construction of $L_M$ in Section
\ref{basicindex}, and explain the claim that this is the
``universal'' case. The idea is that the spectral theory of the Dirac operator
twisted with $L_M$ contains all information which can be obtained using any
kind of Dirac operator, formulated in \cite{SchickICM} as Conjecture 1.5.

The version of the Schr\"odinger-Lichnerowicz Theorem which will be crucial for the results in this paper is the following. It is based on the Schr\"odinger-Lichnerowicz formula (equation \ref{eq:Schroedinger}).
\begin{theorem}\label{theo:Schroedinger}
  If a compact spin manifold $M$ has positive scalar curvature, then
  the spin Dirac operator twisted by any {flat bundle (where the
    fibres are
    vector spaces or more generally modules over a $C^*$-algebra)} is invertible. In particular, the Dirac operator twisted with the Mishchenko line bundle is invertible in this case.
\end{theorem}

Let us therefore make the following
\begin{definition}
  Denote by $\Inv$ the space of Riemannian metrics such that the Dirac
  operator twisted with the Mishchenko bundle is invertible. We call these the
  metrics with a \emph{universally invertible Dirac operator}.   
\end{definition}

By Theorem \ref{theo:Schroedinger} we have $\mathcal{P}\subset \Inv$.

In general, the invertibility of the untwisted Dirac operator
depends on the chosen spin structure. However, it is a basic fact that
this is not so for the universal case: invertibility of the Dirac operator
twisted with the Mishchenko line bundle $L_M$ is independent of the chosen
spin structure, so that $\Inv$ is unambiguously
defined. The idea is as follows: if $s_1$ and $s_2$ are two
  spin structures on $M$, then there is a (graded) real metric line bundle $L
  \to M$ such that the spinor bundles $S_1$ and $S_2$ are related by $S_2=S_1
  \otimes L$. The line bundle $L$ has a canonical flat
    connection, and when we are in
  addition twisting with the Mishchenko line bundle, $L$ can be absorbed into the
  latter bundle at the expense of applying an automorphism of the real group
  $C^*$-algebra. This implies that the Dirac operator on $S_2$ twisted by the
  Mishchenko line bundle is unitarily equivalent to that on $S_1.$ (See \cite[Section 3]{Nitsche} for details.) The invertibility claim now follows immediately.

As stated above, essentially all the tools known to study
$\mathcal{P}$ actually extend to $\Inv$. Indeed, it is a challenging and open
problem to understand the difference between $\mathcal{P}$ and $\Inv$
better.

For us, however, the goal is to transfer information about the homotopy type
of $\mathcal{P}$ to $\mathcal{N}$ or rather the weakly homotopy equivalent
$\mathcal{P}^\sharp$. For example, we want to show that path components of
$\mathcal{P}$ which
belong to distinct path-components of $\Inv$ remain distinct also in $\mathcal{P}^\sharp$. This
would evidently be true if $\mathcal{P}^\sharp\subset \Inv$, i.e.~if
$\RicFlat\subset \Inv$. However, in general this is not true.

Nonetheless, and this is one of the main results of this paper, those Ricci
flat metrics which do not have a universally invertible Dirac operator are
completely isolated from all the metrics with positive scalar curvature.

\begin{theorem}\label{theo:RicFlat-decompose}
  We have a disjoint union decomposition 
  \begin{equation*}
    \RicFlat = \RicFlat^{INV}\amalg \RicFlat^s,
  \end{equation*}
  where $\RicFlat^{INV}:=\Inv \cap \RicFlat,$ $\RicFlat^s:=\RicFlat \setminus \RicFlat^{INV},$ 
  and both $\RicFlat^{INV}$ and $\RicFlat^s$ consist of a union of path-components of $\RicFlat.$ 
  For all metrics in $\RicFlat^{s}$, the universal covering admits a
  non-trivial parallel spinor, and (in particular) the metric has special
  holonomy. On the other hand, no metric in $\RicFlat^{INV}$ has a non-trivial parallel spinor on its universal cover.
\end{theorem}

\begin{theorem}\label{theo:main-decompose}
  We can write $\mathcal{P}^\sharp = \RicFlat^{s}\amalg
  (\mathcal{P}^\sharp\setminus \RicFlat^{s})$, where the former is a union
  of path-components of $\mathcal{P}^\sharp$ and the latter embeds into
  $\Inv$. In particular, all information about the non-triviality of the homotopy
  type of $\mathcal{P}$ which factors through $\Inv$ (e.g.~about path-components
  of $\mathcal P$ which belong to distinct path-components of
  $\Inv$) extends to $\mathcal{P}^\sharp$.

Because the inclusion $\mathcal{P}^\sharp\hookrightarrow \mathcal{N}$ is a weak
  homotopy equivalence (in particular a bijection on $\pi_0$),
  we have an analogous decomposition into unions of path components
  $$\mathcal{N} = \RicFlat^{s}\amalg
  (\mathcal{N}\setminus \RicFlat^{s}),$$
and all information about the non-triviality of the homotopy
  type of $\mathcal{P}$ which factors through $\Inv$ extends to $\mathcal{N}$.
\end{theorem}

  
Below, we will give a number of specific examples of the principle described
in Theorem \ref{theo:main-decompose}.

It should be noted that the existence of a parallel spinor for some metric
does not exclude the possibility that the manifold admits metrics of positive
scalar curvature. For example, simply-connected Calabi-Yau 3-folds are known
to admit both 
positive scalar curvature metrics as well as Ricci-flat metrics with parallel
spinors. The existence of positive scalar curvature on these
  objects is automatic as the $\alpha$-invariant, which is the only obstruction for
  simply-connected spin manifolds to admit positive scalar curvature by \cite{Stolz}, vanishes in real dimension six.

One should remark that no example of a Ricci flat metric without parallel
spinor on its universal covering is known. This means that in the decomposition of
Theorem \ref{theo:RicFlat-decompose}, the part {$\RicFlat^{INV}$} might well be
empty in all cases. This would, by Theorem \ref{theo:main-decompose}, mean that up to weak
homotopy equivalence, we obtain the space of non-negative scalar
curvature metrics $\mathcal{N}$ from the space of positive scalar curvature
metrics $\mathcal{P}$ by just adding a collection of very special components, consisting of
Ricci flat metrics with special holonomy.

The main tools to prove Theorem \ref{theo:main-decompose} from Theorem
\ref{theo:RicFlat-decompose} rely on the theory of
special holonomy. To apply them, we have to establish a link between universal
non-invertibility and the existence of parallel spinors. This is done via the
construction of harmonic spinors. On a closed manifold, it is
elementary to see that the (untwisted) Dirac operator is non-invertible if and
only if there is a non-trivial harmonic spinor. Moreover, if the metric has
non-negative scalar curvature, the Schr\"odinger-Lichnerowicz formula,
Equation \ref{eq:Schroedinger}, implies
that a harmonic spinor is parallel, which in turn forces the holonomy to be special.

The more subtle problem is dealing with a non-invertible
Mishchenko twisted Dirac operator. In 
general,  this does not imply the existence of a non-trivial kernel, because the spectrum of such an
infinite dimensional operator is not in general discrete.

However, we can make use of the fact that the existence of a Ricci flat metric
implies that $\pi_1(M)$ is virtually abelian. Using this, we will show that at
least for some finite dimensional twist bundle, the twisted Dirac operator has
a kernel, i.e.~there exists a ``twisted harmonic spinor''. This will be
sufficient to establish Theorem \ref{theo:RicFlat-decompose} and will be discussed in Section \ref{basicindex}.

We now turn our attention to applications of the above results. We reiterate that almost all known invariants which detect topology in the space ${\mathcal P}$ factor
through the space $\Inv(M)$.
This means that most existing results about the topology of the (moduli) space of positive scalar curvature metrics can be
generalized to non-negative scalar curvature. 
We now present some concrete
examples.

The Kreck-Stolz $s$-invariant is an important tool for studying
path-connectedness of moduli spaces of positive scalar curvature metrics. This
was developed and first used in \cite{KS}. The $s$-invariant is defined for
spin manifolds $M^{4n-1}$ ($n \ge 2$) with vanishing real Pontrjagin classes
and positive scalar curvature. It is an invariant of the path-component in the space of positive scalar curvature metrics. Moreover, if $H^1(M;{\mathbb Z}_2)=0$ (which means the spin structure on $M$ is uniquely determined by the orientation), and $g$ is a positive scalar curvature metric on $M$, then $|s(M,g)| \in {\mathbb Q}$ is an invariant of the path-component in the {\it moduli} space of positive scalar curvature metrics on $M$ containing $g$.

Using Theorem \ref{theo:main-decompose} we can establish:
\begin{theorem}\label{s-inv}
For a closed spin manifold $(M,g)$ of dimension $4k-1,$ ($k \ge 2$), with positive scalar curvature and vanishing real Pontrjagin classes, the Kreck-Stolz $s$-invariant is an invariant of the path-component of non-negative scalar curvature metrics containing $g$. If in addition $H^1(M;{\mathbb Z}_2)=0,$ $|s|$ is an invariant of the path-component containing $[g]$ in the moduli space of non-negative scalar curvature metrics.
\end{theorem}

From Theorem \ref{s-inv} we immediately obtain the following result, which is the non-negative scalar curvature analogue of \cite{KS} Corollary 2.15:
\begin{cor}\label{infinite}
Given any $M$ as in Theorem \ref{s-inv} with $H^1(M;{\mathbb Z}_2)=0,$ the moduli space of non-negative scalar curvature metrics on $M$ has infinitely many path-components.
\end{cor}
 
Besides the $s$-invariant, one can re-visit other types of results for (moduli)
spaces of positive scalar curvature metrics established using index theory,
and making the required adjustments re-state these as results about
non-negative scalar curvature. For example, one can do this with the theorems
about the higher homotopy groups of the (observer moduli) space of positive
scalar curvature metrics established in \cite{HSS}, as these results
rely on the invertibility of a family of Dirac operators which is governed
by the existence or otherwise of harmonic spinors. As a sample result, extending \cite[ Theorem 1.1]{HSS} and \cite[Theorem A]{BERW}, we have

\begin{theorem}\label{theo:HSS}
Given $k \in {\mathbb N} \cup \{0\},$ there is an $N(k) \in {\mathbb N}$ such
that for each $n \ge N(k)$ and each closed spin manifold $M^{4n-k-1}$
admitting a metric $g_0$ with positive scalar curvature, the homotopy group
$\pi_k(\mathcal{N},g_0)$, where $\mathcal{N}$ denotes the space of
non-negative scalar curvature metrics on $M$, contains elements of infinite
order if $k\ge 1,$ and infinitely many different elements if $k=0.$ Their
images under the Hurewicz homomorphism in $H_k(\mathcal{N})$ still have
infinite order.

Indeed, using \cite{BERW}, $N(k)$ can be chosen to be equal to
  $6$ for the statement on $\pi_k$.
\end{theorem}

We expect that also the statement on $H_k$ holds with
  $N(k)=6$. However, the corresponding question for $\mathcal{P}$ is not
  treated in \cite{BERW}.

In precisely the same way, one can generalize to the space $\mathcal{N}$ the classic
results of Hitchin on the non-triviality of $\pi_0(\mathcal{P})$ and
$\pi_1(\mathcal{P})$ for spin manifolds in dimensions 0 and 1, respectively 0
and 7 modulo 8. See \cite{Hi} for the full details, or for a synopsis
explaining the dependence of these results on the invertibility of the Dirac
operator, see IV.7 of \cite{LM}. The same can also be said for the more
recent results of Crowley-Schick (\cite{CS}), Crowley-Schick-Steimle
(\cite{CSS}), Botvinnik-Ebert-Randal-Williams \cite{BERW} and Ebert-Randal-Williams \cite{ERW}, as the underlying analytic facts are precisely the same as in
Hitchin's work.

We also use Theorem \ref{theo:main-decompose} to derive some new examples
involving Ricci non-negative metrics. We remark that the following theorem
presents merely one set of examples among many that are possible. Details of
the Bott manifold $B^8$ appearing in this theorem are given in section
\ref{index}.

\begin{theorem}\label{K3}
If $K^4$ denotes the K3 surface, $B^8$ the Bott manifold, and $\Sigma^{4n-1}$ is any homotopy $(4n-1)$-sphere ($n \ge 2$) which bounds a parallelisable manifold, then both $\Sigma \times K^4$ and $\Sigma \times B^8$ have infinitely many path-components of non-negative Ricci curvature metrics. 
\end{theorem}

As far as the authors are aware, Theorem \ref{K3} is the only result to date concerning the topology of the space of Ricci non-negative metrics in the simply-connected case. It should be noted that we cannot use Theorem \ref{s-inv} to establish these examples as the real Pontrjagin classes are not all zero, and so the $s$-invariant is not defined. The important thing here is that although the manifolds above are known to admit metrics which have both positive scalar and non-negative Ricci curvature, none are known to admit metrics with strictly positive Ricci curvature. There are no known obstructions to positive Ricci curvature for these manifolds: besides admitting positive scalar curvature, they also have finite fundamental group and thus comply with Myers' Theorem.

Since the initial version of this paper was made available, other results
concerning the topology of the moduli space of Ricci non-negative 
metrics have 
appeared, see \cite{TWi}. These results rely on the fundamental group being non-trivial, in contrast to Theorem \ref{K3}. Specifically, \cite{TWi} 
contains examples of manifolds for which the moduli space of Ricci
non-negative metrics has infinitely many path components in both the closed
case (in all dimensions $\ge 7$) as well as in the complete non-compact case (in all
dimensions $\ge 8$). It is also established there that the higher homotopy and rational cohomology groups of the moduli space can be non-trivial in certain cases.

This paper is laid out as follows. In Section \ref{sec:geom_proofs} we collect the geometric
results and prove Theorems \ref{maintheorem-a} and
\ref{theo:main-decompose}. In Section \ref{basicindex} we recall the basic
constructions of higher index theory of (twisted) Dirac operators and harmonic
spinors and prove Theorem \ref{theo:RicFlat-decompose}. In Section \ref{index}
we prove the concrete applications of index theory to spaces of metrics with
non-negative scalar curvature and non-negative Ricci curvature.

This paper grew out of a paper with the same title by the second-named author,
and in relation to this he would like to express his deep gratitude to Bernd
Ammann for his interest and extensive correspondence which considerably
enhanced the paper. Thanks also go to Anand Dessai, Wilderich Tuschmann,
Guofang Wei, Hartmut Weiss and Mark Walsh for their comments. Finally, we
thank an anonymous referee for many helpful comments, in particular for
pointing out a wrong argument in the proof of Proposition \ref{invert} and
providing Proposition \ref{prop:abstract_spec_dect} with its proof as a remedy.


\section{Proofs of the geometric results}\label{sec:geom_proofs}

We want to start with the proofs of our ``geometric'' results, which are
actually independent of the higher index theory discussed at the end of the
introduction.

Our first result is Theorem \ref{maintheorem-a}, which is a rather direct
consequence of the powerful machinery of the Ricci flow. Versions of Theorem
\ref{maintheorem-a} are certainly known to the experts.
As a preliminary, we consider the effects of the Ricci flow on metrics with non-negative scalar curvature.

\begin{lemma}{\cite[2.18]{Br}}\label{flow}  If $M$ is a closed manifold and
  $g_0$ is a metric on $M$ with non-negative scalar curvature, consider the
  Ricci flow $g(t)$ with $g(0)=g_0.$ Suppose that the flow exists for all $t
  \in [0,T].$ Then $g(t)$ has non-negative scalar curvature for all $t \in
  [0,T].$ Moreover, $g(t)$ has positive scalar curvature for all $t \in (0,T]$
  unless $g_0$ is Ricci flat, in which case $g(t)=g_0$ for all $t \in [0,T]$. 
\end{lemma}

\noindent{\bf Proof of Theorem \ref{maintheorem-a}}.
Let $f\colon (D^n,S^{n-1})\to (\mathcal{N},\mathcal{P}^\sharp)$ be
continuous. By \cite[Chapter II, Lemmas (3.1) and (3.2)]{Whitehead} (in conjuction with
\cite[Chapter IV, Section 7]{Whitehead}) we have to find a homotopy
\begin{equation*}
F\colon
(D^n\times [0,T],S^{n-1}\times [0,T])\to
(\mathcal{N},\mathcal{P}^\sharp)\quad\text{
such that } F(D^n\times\{T\})\subset \mathcal{P}^\sharp.
\end{equation*}

Due to the results on the short time existence of the Ricci
flow and the compactness of $D^n$, there is indeed $T>0$ such that the Ricci flow 
defines {a map $F\colon D^n\times [0,T]\to
\mathcal{N}$ with $F|_{D^n\times \{0\}}=f$. By \cite[Theorem A]{BGI}, the Ricci flow depends continuously on the initial data, and thus the map $F$ is continuous.} By Lemma \ref{flow}, $F(D^n\times
(0,T])\subset \mathcal{P}^\sharp$. This means that all the required conditions for the above homotopy are
satisfied. As $f$ was arbitrary, the assertion follows.
\hfill\sq    

We now address our second geometric result, Theorem
\ref{theo:main-decompose}, stating that Ricci flat metrics
which do not have a universally invertible Dirac operator are isolated among metrics
with non-negative scalar curvature.

\begin{lemma}\cite[Satz 2]{Friedrich}\label{Futaki} (See also \cite{Fu}.) If $N$ is a connected Riemannian spin
manifold with a non-zero parallel spinor, then $N$ is Ricci flat.
\end{lemma}

The existence of a parallel spinor on a compact Riemannian spin manifold has
consequences beyond the Ricci flatness of the metric. Indeed, the next result
shows that there cannot be positive scalar curvature metrics arbitrarily
close-by.

\begin{theorem}\label{Dai}
(\cite{DWW}, Theorem 4.2 and subsequent Remark) If $(M,g)$ is a closed
Riemannian spin manifold  with a non-trivial parallel spinor, then there is no
path of metrics $g_t,$ with $g_0=g,$ such that $\mathrm{scal}(g_t)>0$ for all
$t>0$. More generally, there is no path of non-negative scalar curvature metrics
$g_s$ with $g_0=g$ containing a sequence of positive scalar curvature metrics
$g_{s_n},$ where $s_n \xrightarrow{n\to\infty} 0$.
\end{theorem}

The existence of a parallel spinor on a compact Riemannian spin manifold
places restrictions on the holonomy group of that manifold. For a discussion
about these points and detailed references, see for example \cite[Section
1]{AKWW}. Although we will not use holonomy arguments directly, the above
results from \cite{DWW} depend in part on such matters. One might also compare
the results in \cite{Wang}. Holonomy is central to the paper \cite{AKWW}, from
which we will need the following theorem:
\begin{theorem}\label{Cor3}
\cite[Corollary 3]{AKWW} Let $(M,g_0)$ be a closed Riemannian spin manifold
which admits a parallel spinor on its universal cover. If $g_t$, $t \in
[0,T],$ is a smooth family of Ricci-flat metrics on $M$ extending $g_0$, then
the pull-back of $g_t$ to the universal cover admits a parallel spinor for all
$t \in [0,T]$, and the dimension of the space of parallel spinors is
independent of $t$.
\end{theorem}

There is one final result from the literature which we will need, and this is
the basic structure theorem for Ricci-flat metrics (see \cite{CG}, or 4.5 of
\cite{FW}), which also enters crucially in the proof of Theorem \ref{Cor3} above.

\begin{theorem}\label{structure}
(The Ricci-flat structure theorem.) If $(M,g)$ is a closed Ricci-flat
manifold, then there is a finite normal Riemannian covering
$\pi\colon (\bar{M},\bar{g}) \times (T^q,h_{fl}) \to (M,g),$ where
$(\bar{M},\bar{g})$ is a simply-connected Ricci-flat manifold and
$(T^q,h_{fl})$ is the $q$-torus equipped with a flat metric. In particular,
$\pi_1(M)$ contains a free abelian subgroup of finite index.
\end{theorem}

With this preparation at hand, we are now in a position to prove 
Theorem \ref{theo:main-decompose}, assuming Theorem
\ref{theo:RicFlat-decompose}. The essential point is to generalize Theorem
\ref{Dai} from closed manifolds with a parallel spinor to closed manifolds
whose universal covering has a parallel spinor:

\begin{proposition}\label{prop:Dai_plus}
  Let $(M,g_0)$ be a closed Riemannian manifold such that its universal
  covering is spin with a non-zero parallel spinor. Let $ (g_t, 0\le t\le T)$ be a
  continuous path of metrics with $g_t\in \mathcal{P}^\sharp$ starting at $g_0$. Then $g_t\in
  \RicFlat$ for all $t\in [0,T]$.
\end{proposition}

The following result will be used in the proof of Proposition \ref{prop:Dai_plus}.
\begin{lemma}\label{new}
If $(M,g_0)$ is a closed Riemannian manifold such that its universal
  covering is spin with a non-zero parallel spinor, then there exists a finite Riemannian covering $(\bar{M},\bar{g})$ which has a parallel spinor.
\end{lemma}

\begin{proof}
  By Lemma \ref{Futaki}, the existence of a non-zero parallel spinor on the 
  universal covering of $(M,g_0)$ means that the universal cover is Ricci flat, 
  from which it follows that $(M,g_0)$ is also Ricci flat. By Theorem \ref{structure}, some
  finite covering $(\bar M,\bar g)$ of $(M,g_0)$ is a Riemannian product
  $(N,h_N)\times (T^q,h_{fl})$ with simply connected $N$, so the universal covering of $(M,g_0)$ is the
  Riemannian product $(N,h_N)\times (\mathbb{R}^q,h_{fl})$. The existence of a 
  non-zero parallel spinor on a Riemannian product is equivalent to the existence of 
  a parallel spinor on each factor individually, compare
  e.g.~\cite[Theorem 2.5]{L}. In particular, $(N,h_H)$ admits a
  parallel spinor. With a suitable spin structure, $(T^q,h_{fl})$ also has a
  parallel spinor, and we conclude that the closed manifold $(\bar M,\bar g)$
  admits a parallel spinor (with a suitable spin structure).
\end{proof}

\begin{proof}[Proof of Proposition \ref{prop:Dai_plus}] Let $t_1\in [0,T]$ be maximal such that $g_t\in \RicFlat$ for all $t\in
  [0,t_1]$. This exists because $\RicFlat$ is closed. Combining Theorem \ref{Cor3} 
  with the arguments of the above paragraph, we see that $(\bar{M},\bar{g}_{t_1})$ 
  has a parallel spinor (for a suitable spin structure). If $t_1<T$ we could now
  directly apply Theorem \ref{Dai} to the path $(\bar g_t, t_1\le t\le T)$ of
  non-negative scalar curvature metrics lifted to $\bar M$, to deduce that $\bar
  g_t\notin\mathcal{P}$ for $t$ close to $t_1$, $t>t_1$. Therefore $\bar g_t\in
  \RicFlat$, and hence $g_t\in \RicFlat$ for such $t$. This is a
  contradiction to the maximality of $t_1$, so $t_1=T$, and the claim is proved.
\end{proof}

\begin{proof}[Proof of Theorem \ref{theo:main-decompose}] (Assuming Theorem \ref{theo:RicFlat-decompose}.)
  By Theorem \ref{theo:RicFlat-decompose}, every metric in $\RicFlat^s$ is
  such that its universal covering admits a parallel spinor.
  Therefore, by Proposition \ref{prop:Dai_plus}, a path in $\mathcal{P}^\sharp$ which
  starts in $\RicFlat^s$ must remain in $\RicFlat$. But then Theorem
  \ref{Cor3} implies that each metric in the path admits a parallel spinor on
  its universal covering. It now follows from Theorem
  \ref{theo:RicFlat-decompose} that the path remains in $\RicFlat^s$,
  i.e.~$\RicFlat^s$ is a union of path components of $\mathcal{P}^\sharp$.

  The decomposition in Theorem \ref{theo:RicFlat-decompose} shows that $\RicFlat\setminus\RicFlat^{s}\subset \Inv$, and
  therefore by the Schr\"odinger-Lichnerowicz Theorem \ref{theo:Schroedinger} we also have
  $\mathcal{P}^\sharp\setminus \RicFlat^s\subset \Inv$.
\end{proof}


\section{Twisted Index Theory and Harmonic
  Spinors}\label{basicindex}\label{spinor}

In this section, we review some facts about the index theory of Dirac operators on a
spin manifold $M$, potentially twisted with a flat Hermitian bundle, where this flat 
bundle is allowed to be a Hilbert $A$-module bundle for an auxiliary $C^*$-algebra
$A$. We will then also study the theory of harmonic and parallel spinors in
this context, and prove in particular Theorem \ref{theo:RicFlat-decompose}. 

However, we will only use $A$-module bundles in a very special
situation. The
relevant $C^*$-algebra always is the group $C^*$-algebra $C^*\pi$ of the
fundamental group $\pi=\pi_1(M)$ of a Ricci flat manifold $M$. By the structure Theorem
\ref{structure}, $\pi$ then contains a free abelian subgroup of finite index,
and in particular is amenable, so there is only one group $C^*$-algebra: $C^*_{red}\pi=C^*_{max}\pi=:C^*\pi$.

The relevant flat $C^*\pi$-module bundle is the
`Mishchenko line bundle' over $M.$ This is a bundle whose fibre is a free rank one
module over $C^*\pi$, constructed as follows. Let $\tilde M$ be a universal cover of $M.$ There is a free right action
of $\pi$ on $\tilde M$ and a left action on $C^*\pi$, which allows
us to form the flat $C^*\pi$-line bundle
$$ L_M:=\tilde M
\times_{\pi} C^*\pi \to M.$$
 Despite the terminology, if we choose to view this as a
 complex vector bundle, its dimension is equal to the order of $\pi_1(M)$.

 We note here that a \emph{Hilbert} $A$-module structure on an $A$-module
 generalizes the
 Hermitian structure in the case $A=\mathbb{C}$; it consists of an $A$-valued
 inner product satisfying suitable axioms. The basic concepts about Hermitian
 structures generalize readily, compare \cite{Lance}.

 The Mishchenko line bundle is the ``universal'' flat Hilbert $A$-module bundle in a
 precise sense as follows:
 \begin{proposition}\label{prop:flat_structure}
   Let $E\to M$ be any flat Hermitian bundle, or more generally a Hilbert
   $A$-module bundle for some $C^*$-algebra $A$ with fibre a finitely
   generated projective $A$-module $P$ (a Hermitian bundle in the special case
   $A=\mathbb{C}$ and $P=\mathbb{C}^d$). Such a flat bundle corresponds to a
   (holonomy) representation
   $\rho\colon \pi\to U_A(P)$. In the special case of a Hermitian bundle this
   is a unitary representation $\rho\colon \pi\to U(d)$.

  By the universal property of the (maximal) group $C^*$-algebra, this
  representation extends to a $C^*$-algebra homomorphism $\rho\colon
  C^*\pi\to \End_A(P)$, making $P$ a $C^*\pi$-$A$-bimodule (in particular, a
  $C^*\pi$-left module).
  The flat bundle $E$ is then obtained as an associated bundle from the
  Mishchenko line bundle by fibrewise tensor product:
  \begin{equation}\label{eq:def_E}
    E = L_M\otimes_{C^*\pi}P.
  \end{equation}

 If $M$ has a spin structure then the Dirac operator $\mathcal D_E$ twisted by $E$,
 acting as an unbounded operator on the Hilbert $A$-module of $L^2$-sections
 of the spinor bundle twisted by $E$, is
 obtained from the Mishchenko twisted Dirac operator $\mathcal D_{L_M}$ as
 follows: one tensors
 its domain $C^*\pi$-module over $C^*\pi$ with $P$ (and
   completes appropriately), and one tensors the operator
 with the identity,
 \begin{equation}\label{eq:D_E_formula}
 \mathcal D_E=\mathcal D_{L_M}\otimes_{C^*\pi}1_P.
\end{equation}
\end{proposition}
\begin{proof}
  All of this follows directly from the definitions. For \eqref{eq:def_E}
  observe that $E=\tilde M\times_\pi P$. Moreover, $P=C^*\pi\otimes_{C^*\pi}
  P$, so that finally $E=\tilde M\times_\pi C^*\pi\otimes_{C^*\pi}P=
  L_M\otimes_{C^*\pi}P$. Tracing the identifications, this holds with $\pi$
  and $C^*\pi$ both acting on $P$ via $\rho$.

  The statement about the Dirac operators follows again directly from the
  definitions as unbounded Hilbert $A$-module operators, compare \cite{Lance}
  and \cite{Skandalis}.
\end{proof}

This can be used to show that invertibility of the Mishchenko-twisted Dirac operator implies invertibility for all Dirac operators twisted with flat bundles.

\begin{theorem}\label{theo:invertibility_inheritance}
  Let $M$ be a connected spin manifold, $A$ a $C^*$-algebra, $L\to M$ a flat
  bundle with fibers finitely generated projective $A$-modules, with typical
  fiber the $A$-module $P$. This corresponds to a (holonomy) representation
  $\rho\colon\pi\to U_A(P)$. As in Proposition \ref{prop:flat_structure},
  write $L$ as a bundle associated to the Mishchenko bundle $L_M$,
  $L=L_M\otimes_\rho P$.

  The spectrum of the $L$-twisted Dirac operator $\mathcal D_{L}$ is contained
  in the spectrum of the Mishchenko-twisted Dirac operator $\mathcal
  D_{L_M}$. In particular, if $\mathcal D_{L_M}$ is invertible, i.e.~$0$ is
  not in its spectrum, the same is true for $\mathcal D_{L}$.

  If the $C^*$-algebra homomorphism $\rho\colon C^*\pi\to \End_A(P)$ is
  injective, the spectra of $\mathcal D_L$ and $\mathcal D_{L_M}$ even coincide.
\end{theorem}
\begin{proof}
  By Proposition \ref{prop:flat_structure}, $\mathcal D_{L} = \mathcal D_{L_M}\otimes_\rho1_P$.
  The statement about the spectra therefore is a direct 
  consequence of the corresponding general and abstract
  result for spectra of unbounded operators on Hilbert $A$-modules as
  presented in \cite[14.25]{Skandalis}.
\end{proof}

We now turn to the discussion and application of harmonic spinors. By
definition, a harmonic spinor is a section of the spinor bundle belonging to
the kernel of the Dirac operator. Similarly, for a finite dimensional
flat Hermitian bundle $E$, we define an $E$-twisted harmonic spinor as an
element in the kernel of $\mathcal D_E$.

It is a standard fact in the theory of elliptic self-adjoint  operators that, in this
situation, $\mathcal D$ and $\mathcal D_E$ are invertible if and only if there is no non-trivial
(twisted) harmonic spinor.

Note that the situation is more complicated for the Mishchenko-twisted Dirac operator
$\mathcal D_{L_M}$. Typically, if $\pi$ is infinite, even if $0$ is in the spectrum of
$\mathcal D_{L_M}$, its kernel will be trivial due to the presence of a continuous
spectrum in this situation.

For us, harmonic spinors are important because they give rise to
parallel spinors, which we need for our special holonomy considerations. We
first observe that twisted parallel spinors suffice to guarantee the existence
of a regular parallel spinor on the universal covering.

\begin{proposition}\label{untwist} Suppose that a closed Riemannian spin
  manifold $M$ admits a non-zero parallel twisted spinor for some finite dimensional twisting bundle. Then the universal
  cover equipped with the pull-back metric admits a regular non-zero parallel
  spinor.
\end{proposition}
\begin{proof}
Let $S$ denote the spinor bundle on $M$, let $E\to M$ be a flat Hermitian bundle with
corresponding (holonomy) representation $\rho\colon \pi\to U(d)$, and $S
\otimes E$ the twisted spinor bundle.  Suppose that $\nabla\sigma \equiv 0$
for some $\sigma \in \Gamma(S \otimes E)$.

We first observe that the pull-back $\tilde E$ of $E$ to the universal
cover $\tilde{M}$ is a trivial bundle with trivial flat connection: it is
associated to the holonomy representation $\{1\}=\pi_1(\tilde M)\to
\pi_1(M)\xrightarrow{\rho}U(d)$, which is obviously trivial. Consequently, as
a flat bundle $\tilde E\cong \tilde M\times \mathbb{C}^d$.

The pull-back $\tilde\sigma$ of $\sigma$ to $\tilde M$ is a parallel section
of $\tilde S\otimes \tilde E$ with respect to the pull-back
connection. This is the usual twisted spinor connection because the covering
projection is a local isometry, locally preserving all structures. Using the
identification $\tilde E=\tilde M\times \mathbb{C}^d$ (as flat bundles), we
identify $\tilde S\otimes \tilde E$ with $(\tilde S)^d$ (as bundle with
connection), and
$\tilde \sigma$ can be identified with a vector of $d$ parallel spinors on
$\tilde M$. Because $\sigma$ and therefore $\tilde \sigma$ is non-trivial, at
least one of these components is non-trivial, providing a regular non-zero
parallel spinor on $\tilde M$.

\end{proof}

The next lemma is more or less standard. It is a key result for the proof of
Theorem \ref{theo:Schroedinger}.
\begin{lemma}\label{parallel}
Let $(M,g)$ be a closed connected spin manifold with non-negative scalar curvature. Let
$E\to M$ be a flat finite dimensional Hermitian bundle and assume that there
is a non-trivial $E$-twisted harmonic spinor. Then $g$ is Ricci flat and every
twisted harmonic spinor is parallel.
\end{lemma}

\begin{proof}
The main argument needed here is well-known, see for example
\cite[II.8.10,II.8.17-II.8.18]{LM}. It begins with the
Schr\"odinger-Lichnerowicz formula
\begin{equation}\label{eq:Schroedinger}
{\mathcal
  D_E}^2=\nabla^*\nabla+\frac{1}{4}\mathrm{scal},
\end{equation}
where ${\mathcal D_E}$ is the
twisted Dirac operator and $\nabla^*\nabla$ is the connection Laplacian on
spinors twisted by the flat bundle $E$ with its flat connection. Because the
connection of $E$ is flat, there is no additional term on the right hand
side. Given any non-trivial $E$-twisted harmonic spinor $\sigma$, integrating
over $M$  gives the following equation:
$$\int_M
\frac{\mathrm{scal}\cdot|\sigma|^2}{4}+|\nabla \sigma|^2=0,$$ where the form of the
second term uses the definition of the connection
Laplacian $\nabla^*\nabla$. Thus in the context of non-negative scalar
curvature, we see that $|\nabla\sigma| \equiv 0$ on $M$ and, and thus
$\sigma$ is a non-trivial parallel $E$-twisted
spinor. By Proposition \ref{untwist}, there is a non-zero parallel spinor on
the universal covering of $M$. As the existence of a parallel spinor forces the
metric to be Ricci flat by Lemma \ref{Futaki}, the universal covering of $M$, and
therefore $M$ itself, are both Ricci flat.
\end{proof}  

The final preparational result provides twisted harmonic spinors if the
metric does not have a universally invertible Dirac
operator, but only in the case of our very special
fundamental group. This is a partial converse to Theorem
\ref{theo:invertibility_inheritance}, and is probably well known.

\begin{proposition}\label{invert} Let $M$ be a closed Riemannian spin manifold such that
  its fundamental group $\pi$ has a free abelian subgroup of finite index.
  If for every finite dimensional flat Hermitian bundle $E$ the twisted Dirac
  operator ${\mathcal D}_E$ is invertible, then the metric
    has a universally
  invertible Dirac operator. Equivalently, if the metric
  does not have a universally invertible Dirac operator then it
  admits a non-zero twisted harmonic spinor.
\end{proposition}

The proof relies on the following detection principle for the spectrum of Hilbert
$C^*\pi$-module operators. This is our
  interpretation of the classical Floquet-Bloch theory. It was provided to us,
  together with its proof, by the anynomous referee, and we are grateful for
  this help.




\begin{proposition}\label{prop:abstract_spec_dect}
 Assume that the group $\pi$ contains the subgroup $\integers^n$ with finite
  index $d$. Let $a$ be a possibly unbounded self-adjoint Hilbert
  $C^*\pi$-module 
  operator on the countably generated Hilbert $C^*\pi$-module $E$, such that
  its bounded transform $T:= a(a^2+1)^{-1/2}$ satisfies the property that
  $S:=T^2-1$ is
  compact in the sense of Hilbert $C^*\pi$-module morphisms.

If $0$ is in the spectrum of $a$, then $0$ is
  already in the spectrum of $a\tensor_\rho 1$ for at least one
  representation $\rho\colon \pi\to U(d)$. 
\end{proposition}
\begin{proof}
By the spectral mapping theorem, $0$ is in the spectrum of $a$ if and only
  if $0$ is in the spectrum of $T$ if and only if $-1$ is in the spectrum of $S$, so we study $S$ instead of $a$.
  
Next, we normalize the Hilbert module by taking the direct sum of $S$ with the zero operator on $l^2(C^*(\pi))$. By
  Kasparov's stability theorem, we then
      may assume that $E = l^2(C^*(\pi))$.
  
Assume initially that $\pi=\integers^n$, with Fourier transform
  isomorphism $C^*\integers^n\to 
  C(T^n)$. Recall that $C^*\integers^n$ is the $C^*$-algebra of bounded operators on
  $l^2(\integers^n)$ generated by convolution with $z_i$,
  where $z_1,\dots,z_n$ are generators of the infinite cyclic
  summands. The Fourier transform isomorphism $l^2(\integers^n)\cong
  L^2(T^n)$ just reinterprets $z_i$ as a variable of the factor
  $S_i^1\subset
  \complexs$ of the torus $T^n$. Under this identification, convolution
  with $z_i$ becomes multiplication by $z_i$, which is now a continuous
  function on $T^n$. In this way, $C^*\integers^n$ is identified with a
  $C^*$-subalgebra of $C(T^n)$, and by the Weierstra\ss\ approximation theorem
  is indeed all of $C(T^n)$.

   Next, when passing to  Hilbert $C(T^n)$-modules we have the isomorphisms
  $\oplus_{k\in\mathbb{N}}
  C(T^n)= l^2(C(T^n))\iso C(T^n, l^2)$, with the $C(T^n)$-valued inner-product
  defined pointwise. The crucial fact now is that the $C^*$-algebra of compact
  $C(T^n)$-Hilbert module operators is identified with 
  $C(T^n, K(l^2))$, where $K(l^2)$ is the algebra of compact operators on the
  Hilbert space $l^2$ with the norm topology. We thank the referee for pointing
  out that the corresponding statement is not true for the bounded operators,
  when using the norm topology on $B(l^2)$.

  For a norm continuous function valued in compact operators $S\in C(T^n,K(l^2))$, it is
  clear that $S-\lambda$ is invertible if and only if for each $\rho\in T^n$
  the operator $S(\rho)-\lambda=S\tensor_\rho 1-\lambda $ is invertible. This
uses the fact that the subset of 
  invertible operators on $l^2$ is open in $B(l^2)$, and we use the
  interpretation of $\rho\in T^n$ as evaluation {homomorphism $\rho\colon
  C(T^n)\to\complexs$}.

  Now we pass to the general situation, i.e.~$\integers^n$ is a subgroup of finite index of $\pi$.
 Choose a set $\{g_1,\dots,g_d\}$ of right coset representatives for
  $\integers^n$ in $\pi$. We obtain the Fourier isomorphism
  \begin{equation*}
  l^2(\pi)=\oplus_{j=1}^d g_j l^2(\integers^n) \iso \oplus g_j
  L^2(T^n)=L^2(T^n, \oplus_{j=1}^d g_j\complexs),
\end{equation*}
with $L^2(T^n,\oplus g_j\complexs)$ the space of
  $\complexs^d$-valued $L^2$-functions on $T^n$. 

  Left multiplication by an element $g\in \pi$ permutes the right cosets and
  maps $g_j$ to $g_{\alpha(j)}v_j$ with $v_j\in\integers^n$, ($v_j$ and the
  permutation $\alpha$ depend on $g$). Under our Fourier transform
  isomorphism, this operator becomes the 
  operator which multiplies the $j$-th component with the Fourier polynomial
  $v_j\in C(T^n)$, and then applies pointwise the permutation matrix
  $\alpha$. In particular, the closure, $C^*\pi$, is identified with a
  sub-$C^*$-algebra of the matrix-valued continuous functions $C(T^n,
  M_d(\complexs))$, which we interpret as the $C^*$-algebra of
  endomorphisms of the Hilbert $C(T^n)$-module $C(T^n)^d$.

 The inclusion $C^*\pi\hookrightarrow \End_{C(T^n)}(C(T^n)^d)$
   allows us to induce the Hilbert $C^*\pi$-module $l^2(\pi)$ up to the Hilbert
   $C(T^n)$-module $l^2(C^*\pi)\tensor_{C^*\pi} C(T^n)^d\iso l^2(C(T^n)^d)$. This gives rise to
   the embedding
   $\End_{C^*\pi}(l^2(C^*\pi))\hookrightarrow \End_{C(T^n)}(l^2(C(T^n)^d)); S\mapsto
   S\tensor 1_{C(T^n)^d}$ which maps compact elements to compact elements. By
   \cite[14.25]{Skandalis}, used already in the proof of Theorem 
  \ref{theo:invertibility_inheritance}, under this embedding the spectrum is
  unchanged.  Using the special case
  of $\integers^n$ we already established, the spectrum is then detected by looking at the
  induced operators $(S\tensor 1_{C(T^n)^d})\tensor_\rho 1_\complexs$
    for the evaluation homomorphisms
  $\rho\colon T^n\to\complexs$, because $l^2(C(T^n)^d)\iso l^2(C(T^n))$.

Composed with the embedding $C^*\pi\hookrightarrow
  C(T^n,M_d(\complexs))=\End_{C(T^n)}(C(T^n)^d)$, such an evaluation
  homomorphism becomes the
  homomorphism associated to the representation $R\colon \pi\to
  M_d(\complexs)= \End(\oplus g_j\complexs)$ induced up from the 
  irreducible representation of $\integers^n$ corresponding to $\rho$. This is
  true because, by definition, in this induced 
  representation $g\in\pi$ maps the basis element $g_j$ of $\oplus
  g_j\complexs$ to  $g_{\alpha(j)}\rho(v_j)$, if $gg_j=g_{\alpha(j)}v_j$ as
  above. Consequently, the spectrum of $S$ is detected by looking at the
spectrum of the operators $S\tensor 1_{C(T^n)^d}\tensor_\rho 1_\complexs=
  S\tensor_R1_{\complexs^d}$ for the induced representations $R\colon \pi\to
  M_d(\complexs)$.

To wrap up: if $0$ is in the spectrum of $a$, then $-1$ is
    in the spectrum of $S$, i.e.~$S+1$ is not invertible. This implies that
    $(S+1)\tensor_\rho 1$ is not invertible for some representation
    $\rho\colon \pi\to U(d)$. Because induction is compatible with functional
    calculus, this implies finally by the spectral mapping theorem that $0$ is
    in the spectrum of $a\tensor_\rho 1$.
   \end{proof}

\begin{proof}[Proof of Proposition \ref{invert}]
  We now deal with the unbounded self-adjoint operator $\mathcal{D}_{L_M}$ on
  the Hilbert 
  $C^*\pi$-module of sections of the Mishchenko-twisted spinor bundle. By
  definition, the metric does not have a universally invertible Dirac operator if
  $0$ is in the spectrum of $\mathcal{D}_{L_M}$.

  We can apply Proposition
  \ref{prop:abstract_spec_dect} to $\mathcal{D}_{L_M}$ because the bounded transform
  $T:=\mathcal{D}_{L_M} (\mathcal{D}^2_{L_M}+1)^{-1/2}$ is known by elliptic theory to be a
  bounded self-adjoint Hilbert $C^*\pi$-module operator such that $T^2-1$ is
  compact in the Hilbert $C^*\pi$-module sense.

  Therefore,
  Proposition \ref{prop:abstract_spec_dect} together with Theorem  
  \ref{theo:invertibility_inheritance} on the identification of
  $\mathcal{D}_{L_M}\tensor_\rho id$ with
   $\mathcal{D}_\rho=\mathcal{D}_L$, the Dirac operator twisted with the
  flat bundle $L$ associated to the representation $\rho\colon \pi\to U(d)$,
  imply 
that if $0$ is in the spectrum of $\mathcal{D}_{L_M}$ then $\mathcal{D}_L$ is not invertible for
  at least one finite dimensional flat bundle $L$.
\end{proof}

\begin{proof}[Proof of Theorem \ref{theo:RicFlat-decompose}] 
  We begin by arguing that all metrics in $\RicFlat^s$ admit a non-trivial 
  parallel spinor on the universal cover. By Proposition \ref{invert}, if $g\in\RicFlat^s$ then it
  admits a non-zero twisted harmonic spinor. By Lemma \ref{parallel} this
  twisted spinor is parallel, which implies by Proposition \ref{untwist} the
  desired (regular) parallel spinor on the universal covering. It is a
  standard result that the existence of a non-zero parallel spinor forces the
  holonomy group to be special, compare e.g.~\cite{AKWW}.

  Next, we argue that the existence of a parallel spinor on the universal
  covering implies that the metric does not have a universally
    invertible Dirac operator. Observe that
  by the Schr\"odinger-Lichnerowicz formula \eqref{eq:Schroedinger}, because
  the Ricci and therefore scalar curvature now vanish identically, every
  (twisted) parallel 
  spinor which is square integrable lies in the kernel of the Dirac operator,
  and this applies in particular to every (twisted) parallel spinor on a finite covering of $M$.

By Lemma \ref{new}, the parallel
  spinor on the universal covering gives rise to a parallel spinor on a
  suitable finite covering. A priori, this is for a spin structure different
  from the one pulled back from $M$. But the spinor bundle for this a priori
  different spin structure equals the spinor bundle for the pull-back spin
  structure twisted with an appropriate flat line bundle, by \cite[Section
  3]{Nitsche}. Therefore a non-zero parallel spinor on the universal covering
  produces a twisted parallel spinor on a finite covering which by Theorem
  \ref{theo:invertibility_inheritance} implies that the metric
  does not have a universally invertible Dirac operator.

  Finally, given the existence, respectively non-existence, of non-trivial 
  parallel spinors on the universal cover for metrics in $\RicFlat^s,$ respectively 
  $\RicFlat^{INV},$ we note that by Theorem \ref{Cor3} there can be no path 
  within $\RicFlat$ linking
  $\RicFlat^{INV}$ and $\RicFlat^s.$ Hence $\RicFlat^{INV}$ and $\RicFlat^s$ 
  must each be a union of path-components of $\RicFlat.$  
\end{proof}


\section{Applications via Index Theory}\label{index}

In this section we show how to reduce the proof of Theorems \ref{s-inv} and \ref{K3}
to their counterparts for the space $\mathcal{P}$ of positive scalar curvature
metrics, using Theorems \ref{maintheorem-a} and \ref{theo:main-decompose}.

\begin{proof}[Proof of Theorem \ref{s-inv}] 
We have to show that
  for a given path-component 
$C_{\mathcal{N}}$  of
$\mathcal{N}$ which contains several path components $C_1,\dots,C_k$ of
$\mathcal{P}$, the Kreck-Stolz invariants of $C_1,\dots,C_k$ coincide.

Let $g_t$, $t\in [0,1]$ be a path in $C_{\mathcal{N}}$ joining say $C_1$ and
$C_2$. By Theorem \ref{maintheorem-a} we can deform this path slightly to
obtain a path $\tilde g_t$ now in $\mathcal{P}^\sharp$ joining $C_1$ with
$C_2$, as $C_1,C_2$ are open in $\mathcal{P}^\sharp$. Concretely, we obtain
$\tilde g_t$ by applying the Ricci flow for a short time, starting with the
metrics in $g_t$. The path
$\tilde g_t$ lies in $\mathcal{P}^\sharp\setminus \RicFlat^s$ because its path
component contains metrics of positive scalar curvature, in contrast to
$\RicFlat^s$. The arguments in \cite{KS} can now be applied to the path
$\tilde g_t$, using solely that it runs through the space $\Inv$ of metrics
with invertible Dirac operator. This gives the desired invariance properties
of $s$.
\end{proof}

We remark that it is not difficult to show that the (untwisted) Dirac operator is invertible for any metric in $C_{\mathcal{N}},$ and hence the Kreck-Stolz invariant (which makes sense irrespective of the curvature) takes the same value for all metrics in this path-component, not just those with positive scalar curvature.

Turning our attention to examples, we consider two families of products, one
involving a K3 surface $K^4,$ and the other involving a Bott manifold $B^8$ as
a factor. Recall that, as a smooth manifold, $K^4$ can be defined
by $$K^4:=\{(z_0,z_1,z_2,z_3)\,|\,z_0^4+z_1^4+z_2^4+z_3^4=0\} \subset {\mathbb
  C}P^3.$$ 
This is known to support a Ricci flat metric, see for example \cite[page
128]{Be}, though since $\hat{A}(K^4)=-2$ there is no metric of positive scalar
curvature.

A Bott manifold is a closed simply-connected 8-dimensional
    spin manifold
  $B^8$ with $\hat{A}(B^8)=1$, which therefore does not admit a metric of
  positive scalar curvature. We consider here an example constructed by
  D.~Joyce in \cite{J} which has $\mathrm{Spin}(7)$-holonomy, and thus admits
  a Ricci flat metric.

We also consider the set of homotopy spheres which bound parallelisable
manifolds in dimensions $4n-1$, ($n\ge 2$). Although finite for each $n$, the
order of this family grows more than exponentially with dimension. The moduli
space of positive Ricci curvature metrics for each of these spheres was shown
to have infinitely many path-components in \cite{Wr}. This result was
established by exhibiting an infinite family of Ricci positive metrics on each
sphere, and showing that these metrics can be distinguished by their
$s$-invariants.

\begin{proof}[Proof of Theorem \ref{K3}.] 
It suffices to consider $\Sigma^{4n-1} \times K^4$ for some choice of homotopy sphere $\Sigma^{4n-1}$ bounding a parallelisable manifold, as the argument for $\Sigma^{4n-1} \times B^8$ is identical.

In \cite{Wr} it was shown that we can find a sequence of Ricci positive metrics $g_i$ on $\Sigma$ such that $s(\Sigma,g_i) \neq s(\Sigma,g_j)$ whenever $i \neq j$, so $g_i$ and $g_j$ belong to different path-components of the moduli space of positive scalar curvature metrics on $\Sigma$. For each $i$ there is a parallelisable bounding manifold $W_i$ for $\Sigma$ such that $g_i$ extends to a positive scalar curvature metric $\bar{g}_i$ over $W_i$ (product near the boundary), see \cite[Corollary 6.4]{Wr}.

The $W_i$ are constructed by plumbing $D^{2n}$-bundles over $S^{2n}$. If we
consider the oriented union $W_i \cup_\Sigma (-W_j),$ it is established for
example in \cite[page 73]{Ca} that $\hat{A}(W_i \cup_{\sigma} (-W_j))$ is a non-zero
multiple of the difference of signatures
$\mathrm{sig}(W_i)-\mathrm{sig}(W_j).$ As noted in \cite[\S2]{Wr}, for $i
\neq j$ we have $\mathrm{sig}(W_i) \neq \mathrm{sig}(W_j),$ and thus
$\hat{A}(W_i \cup_{\Sigma} (-W_j)) \neq 0.$ As the $\hat{A}$-genus is multiplicative
for products and $\hat{A}(K^4) \neq 0$, we deduce that $$\hat{A}((W_i \times
K^4) \cup_{\Sigma\times K^4} (-W_j \times K^4)) \neq 0.$$

Let $g_K$ denote a Ricci flat metric on $K^4,$ and consider the product
metrics $g_i+g_K.$ These have non-negative Ricci curvature and positive scalar
curvature. By the above, these metrics can be extended to positive scalar
curvature metrics $\bar{g}_i+g_K$ on $W_i \times K^4$. We can
  now use standard results about the Atiyah-Patodi-Singer index, compare
  e.g.~the survey \cite[Section 5]{TW} for background, to obtain
$$\mathrm{ind}(\mathcal{D}^+(W_i \times K^4,
\bar{g}_i+g_K))=\mathrm{ind}(\mathcal{D}^+(W_j \times K^4, \bar{g}_j+g_K))=0.$$
For $i \neq j$ suppose the metrics $g_i+g_K$ and $g_j+g_K$ belong to the same
path-component of non-negative scalar curvature metrics on $\Sigma \times
K^4$, i.e. there is a path $h_t$, $t \in [0,1]$, with $\mathrm{scal}(h_t) \ge
0,$ $h_0=g_i+g_K$ and $h_1=g_j+g_K$. By Theorems \ref{maintheorem-a} and
\ref{theo:main-decompose} (concretely, via application of the Ricci flow) we
can assume that $h_t\in\Inv$, using that $h_0,h_1\in \mathcal{P}\subset
\Inv$ and that $\mathcal{P}$ is open in $\mathcal{N}$.
Let $\bar{h}_t$ be any path of metrics on $W_i \times K^4$, starting with $\bar{g}_i+g_K$, which extend $h_t$ (and take the form of a product near the boundary). 
{By standard index theory arguments, the invertibility of the
  boundary Dirac operator along the path $h_t$ ensures that
   $\mathrm{ind}(\mathcal{D}^+(W_i \times K^4, \bar{h}_t))$, is independent of $t$. Moreover, since the path $\bar{h}_t$ begins with a positive scalar curvature metric, we see that in fact $\mathrm{ind}\mathcal{D}^+(W_i \times K^4, \bar{h}_t)=0$ for all $t$. It then follows from \cite{APS} that} 
\begin{equation*}
\begin{split}
0=&\hbox{ind}(\mathcal{D}^+(W_i \times K^4, \bar{h}_1))-\hbox{ind}(\mathcal{D}^+(W_j \times K^4, \bar{g}_j+g_K)) \\ 
=& \hat{A}((W_i \times K^4) \cup_{\Sigma\times K^4} (-W_j \times K^4)) \\ \neq&0,
\end{split}
\end{equation*}
and we have a contradiction. Thus $g_i+g_K$ and $g_j+g_K$ cannot belong to the same path-component of non-negative scalar curvature metrics, and hence must belong to different path components of Ricci non-negative metrics. 
\end{proof}

As remarked in the introduction, one can replace the homotopy spheres in
Theorem \ref{K3} with other manifolds. For example one could use the infinite
family of 7-dimensional Einstein manifolds $M_{k,l}$ considered in \cite{KS},
which were shown to have infinitely many path-components of Ricci positive
metrics in \cite{KS}, and infinitely many path components of non-negative
sectional curvature metrics in \cite{KPT}.
\medskip


\end{document}